\documentclass[preprint,11pt]{elsarticle}
\usepackage{color,amsmath,amssymb,graphicx,algorithmic,alltt,xy,fullpage,amsthm,wrapfig,setspace,cancel}
\definecolor{refkey}{rgb}{1,0.5,0.5}
\definecolor{labeled}{rgb}{0.5,1,0.5}
\definecolor{Hydroureter}{rgb}{0,0.45,0}
\definecolor{Mistakable}{rgb}{0,0,0.75}
\definecolor{MyDarkRed}{rgb}{0.9,0,0}
\date{}
\input xy
\xyoption{all}
\newtheorem{thrm}{Theorem}
\newtheorem{lem}{Lemma}

\newtheorem{rmk}{Remark}

\newtheorem{prop}{Continuation Result}

\newcommand{\LC}{\kappa}

\newcommand{\X} {{\bf x }}
\newcommand{\x} {{\bf x }}

\newcommand{\Beginproof}{\vspace{0mm} \parindent=0pt
         {\bf Proof.} \hspace{3mm} \parindent=3ex}
\newcommand{\Endproof}{$\Box$ \vspace{5mm}
                        \parindent=3ex}
\newcommand{\psiex}{\psi_{\rm explicit}}
\newcommand{\psiexalpha}{\psi_{\rm explicit, \alpha}}
\newcommand{\Pcr}{P_{\rm cr}}
\newcommand{\rev} { }
\DeclareMathOperator*{\argmax}{arg\,max}

\begin{document}
\begin{frontmatter}
\title{Nonlinear-damping continuation of the nonlinear Schr\"odinger equation -- a numerical study}

\author{G. Fibich}\ead{fibich@tau.ac.il}\author{M. Klein}\ead{morankli@tau.ac.il}
\address{School of Mathematical Sciences, Tel Aviv University, Tel Aviv 69978, Israel}
\begin{abstract}
We study the nonlinear-damping continuation of singular
solutions of the critical and supercritical NLS. Our simulations
suggest that for generic initial conditions that lead to collapse in the undamped NLS, 
the solution of the weakly-damped NLS 
$$
i\psi_t(t,\X)+\Delta\psi+|\psi|^{p-1}\psi+i\delta|\psi|^{q-1}\psi=0,\qquad0<\delta \ll 1,
$$
is highly asymmetric with respect to the singularity time,
and the post-collapse defocusing velocity of the singular core goes
to infinity as the damping coefficient~$\delta$ goes to zero. In the special
case of the minimal-power blowup solutions of the critical NLS,
the continuation is a minimal-power solution with a higher (but
finite) defocusing velocity, whose magnitude increases monotonically
with the nonlinear damping exponent~$q$.
\end{abstract}
\end{frontmatter}
\section{Introduction}
The nonlinear Schr\"odinger equation (NLS)
\begin{equation}\label{eq:DdimensionalNLS}
i\psi_t(t,\X)+\Delta\psi+|\psi|^{p-1}\psi=0,\qquad\psi_0(0,\X)=\psi_0(\X)\in
H^1,
\end{equation}
where~$\X=(x_1,...,x_d)\in\mathbb{R}^d$
and~$\Delta=\partial_{x_1x_1}+\cdot\cdot\cdot\partial_{x_dx_d}$, is
one of the canonical nonlinear equations in physics, arising in
various fields such as nonlinear optics, plasma physics,
Bose-Einstein condensates (BEC), and surface waves. When~$(p-1)d<4$,
the NLS is called subcritical. In that case, all~$H^1$~solutions
exist globally. In contrast, both the critical NLS~$(p-1)d=4$ and
the supercritical NLS~$(p-1)d>4$ admit singular solutions. Since
physical quantities do not become singular, this implies that some
of the terms that were neglected in the derivation of the NLS,
become important near the singularity.

The continuation of NLS solutions beyond the singularity has been an
open question for many years. In 1992, Merle~\cite{Merle-92a}
presented a continuation of the explicit blowup
solutions~$\psiexalpha$ of the critical
NLS,{\rev{see~(\ref{eq:ExplicitBlowupSolutionAlpha})}}, which is
based on slightly reducing the power ($L^2$ norm) of the initial condition.
This continuation has two key properties:
\begin{enumerate}
\item \textbf{Property 1}: The solution is symmetric with respect to the singularity
time~$T_c$.\label{Item:WeakSolProp1}
\item \textbf{Property 2}: After the singularity, the solution can only be determined up to multiplication by a constant
phase term~$e^{i\theta}$.\label{Item:WeakSolProp2}
\end{enumerate}
More recently, Merle, Raphael and
Szeftel~\cite{Merle-Raphael-Szeftel} generalized this continuation
result to Bourgain-Wang singular solutions~\cite{Bourgain-97}. Note,
however, that both the explicit solutions~$\psiexalpha$ and the
Bourgain-Wang solutions are unstable.

In~\cite{Merle-92b}, Merle presented a different continuation, which
is based on the addition of nonlinear saturation. Merle showed
that, generically, as the nonlinear saturation coefficient goes to
zero, the limiting solution beyond~$T_c$ can be decomposed into two
components: A $\delta$-function singular core that extends for $T_c  \le t \le T^0$, and a regular
component elsewhere.

In~\cite{Tao-2009}, Tao proved the global existence and uniqueness
in the semi Strichartz class for solutions of the critical NLS.
Intuitively, these  solutions are formed by solving the equation in
the Strichartz class whenever possible, and deleting any power that
escapes to spatial or frequency infinity when the solution leaves
the Strichartz class. These solutions, however, do not depend
continuously  on the initial conditions, and are thus not a
well-posed class of solutions. Recently, Stinis~\cite{Stinis-2010}
studied numerically the continuation of singular NLS solutions using
the t-model approach.

In~\cite{Fibich-Klein-2010} we analyzed asymptotically and
numerically four potential continuations of singular NLS solutions:~1)~a sub-threshold power continuation, 2)~a
shrinking-hole continuation for ring-type solutions, 3)~a vanishing
nonlinear-damping continuation, and 4)~a complex Ginzburg-Landau
(CGL) continuation. Our main findings were as follows:
\begin{enumerate}
\item The non-uniqueness of the phase of the singular core beyond the singularity (Property~2) is a
universal feature of{\rev{NLS}} continuations.
\item The symmetry with respect to the singularity
time (Property~1) holds if the continuation model is time reversible
and if it leads to a point singularity (i.e., if it defocuses
for~$t>T_c$). Therefore, it is a non-generic feature.
\end{enumerate}
Recently, the post-collapse loss-of-phase phenomena was demonstrated experimentally for intense laser beams propagating in water~\cite{Bonggu-2012}.

In this paper we further study the effect of small nonlinear-damping
in the NLS
\begin{equation}\label{eq:DampedNLS}
i\psi_t(t,\X)+\Delta\psi+|\psi|^{p-1}\psi+i\delta|\psi|^{q-1}\psi=0,\qquad0<\delta\ll1.
\end{equation}
The addition of small nonlinear-damping is physical. Indeed, in
nonlinear optics, experiments suggest that arrest of collapse is
related to plasma formation, and nonlinear damping is
used as phenomenological model for multi-photon absorption by
plasma. In BEC, a quintic nonlinear damping term corresponds to
losses from condensate due to three-body inelastic recombinations~\cite{Markowich-04}.
In addition, the nonlinear-damping term appears in the complex-Ginzburg-Landau (CGL) equation,
which arises in a models of  chemical turbulence,
Poiseuille flow, Rayleigh-B\'ernard convection, Taylor-Couette flow, and 
superconductivity,

In~\cite{Fibich-Klein-2010} we analyzed the continuation of the critical NLS with a vanishing critical nonlinear damping,
i.e., equation~(\ref{eq:DampedNLS}) with~$p=q=1+4/d$. Since the
NLS~(\ref{eq:DampedNLS}) is not time reversible, its solutions are
asymmetric with respect to the time~$T_{\rm arrest}^{(\delta)}$ at which
the collapse is arrested. In particular, in the limit~$\delta_n \rightarrow0+$, the
continuation of~$\psiexalpha(t,r)$
is~$e^{i\theta}\psi^\ast_{\rm{explicit,\kappa\alpha}}(2T_c-t,r)$,
where~$\kappa\approx1.614$. Hence, the defocusing
velocity~$\kappa\alpha$ is higher then the focusing
velocity~$\alpha$. When the initial condition leads to a loglog
collapse in the undamped critical NLS, asymptotic analysis and
numerical simulations suggest that the singular core expands beyond
the singularity at a velocity that goes to infinity
as~$\delta\rightarrow0+$.

The question that we address in this study is whether and how the
results of~\cite{Fibich-Klein-2010} for~$q=p=1+4/d$ will change in
the following cases:
\begin{enumerate}
\item The critical NLS with a supercritical damping exponent~(i.e., $q>p=1+4/d$).
\item The supercritical NLS with~$q\geq p>1+4/d$.
\end{enumerate}

The paper is organized as follows. In Section~\ref{sec:NLS_Review}
we provide a short review of NLS theory. In
Section~\ref{sec:NL_DampingEffect} we review previous rigorous,
asymptotic, and numerical results on the effect of damping in the
NLS. In Section~\ref{sec:CriticalDampingExponent} we show
numerically that in the supercritical NLS, the nonlinear damping
exponent~$q$ has to be strictly higher than the nonlinearity
exponent~$p$, in order to arrest the collapse. This is different from
the critical case, where collapse is arrested for~$q \ge p$. In
Section~\ref{sec:SuperCriticalNLS} we show that solutions of the
supercritical NLS with a small nonlinear damping are asymmetric with
respect to the arrest-of-collapse time~$T_{\rm arrest}^{(\delta)}$, and
that the post-collapse defocusing velocity of the singular core goes to infinity as the
damping coefficient~$\delta$ goes to zero. 
In Section~\ref{sec:CriticalNLS} we obtain similar results
for the critical NLS with generic initial conditions that lead to a
loglog collapse. In the special case of the minimal-power explicit blowup
solution~$\psiexalpha(t,r)$ of the critical NLS, however, the
continuation beyond the singularity is also defined for $q<p$, 
and is given by~$e^{i\theta}\psi^\ast_{\rm{explicit,\kappa(q)\alpha}}(2T_c-t,r)$,
where~$\kappa(q)$ increases monotonically with~$q$. Final remarks are given in Section~\ref{sec:final}.

Overall, the qualitative effect of small nonlinear damping on the collapse  is  the same
in the critical and supercritical  
NLS. One difference is that in the critical case collapse is arrested for~$q \ge p$,
whereas in the supercritical case collapse is only arrested for~$q > p$. Another difference is that
the distance between the damped solution around~$T_{\rm arrest}^{(\delta)}$ and the
asymptotic profile of the undamped NLS is small in the critical case, 
but large in the supercritical case.
Surprisingly, in the latter case, the profile near~$T_{\rm arrest}^{(\delta)}$ appears to be
given by a rescaled supercritical standing wave.

\section{Review of NLS theory}\label{sec:NLS_Review}

 The
NLS~(\ref{eq:DdimensionalNLS}) has two important conservation laws:
{\em {Power}} conservation\footnote{We call the~$L^2$ norm the
{\em{power}}, since in optics it corresponds to the beam's power.}
\[
P(t)\equiv P(0),\qquad P(t)=\int|\psi|^2d\X,
\]
and {\em{Hamiltonian}} conservation
\begin{equation}\label{eq:HamiltonianConservation}
H(t)\equiv H(0),\qquad
H(t)=\int|\nabla\psi|^2d\X-\frac2{p+1}\int|\psi|^{p+1}d\X.
\end{equation}

The NLS~(\ref{eq:DdimensionalNLS}) admits the waveguide
solutions~$\psi=e^{it}R(r)$, where~$r=|\x|$, and~$R$ is the solution
of
\begin{equation}\label{eq:dDim_R_ODE}
R''(r)+\frac{d-1}{r}R'-R+R^p=0,\qquad R'(0)=0,\quad R(\infty)=0.
\end{equation}
When~$d=1$, the solution of~(\ref{eq:dDim_R_ODE}) is unique, and is
given by
\begin{equation}\label{eq:1D_ground_state_solutions}
R_p(x)=\left(\frac{p+1}{2}\right)^{1/(p-1)}\mbox{cosh}^{-2/(p-1)}\left(\frac{p-1}{2}
x\right).
\end{equation}
When~$d\geq2$, equation~(\ref{eq:dDim_R_ODE}) admits an infinite
number of solutions. The solution with the minimal power, which we
denote by~$R^{(0)}$, is unique, and is called the ground state.

\subsection{Critical NLS}

In the critical case~$(p-1)d=4$, equation~(\ref{eq:DdimensionalNLS})
can be rewritten as
\begin{equation}
   \label{eq:dDim_CriticalNLS}
i\psi_t(t,\X)+\Delta\psi+|\psi|^{4/d}\psi=0,\qquad\psi_0(0,\X)=\psi_0(\X)\in
H^1,
\end{equation}
and equation~(\ref{eq:dDim_R_ODE}) can be rewritten as
\begin{equation}\label{eq:dDimCritical_R_ODE}
R''(r)+\frac{d-1}{r}R'-R+R^{4/d+1}=0,\qquad R'(0)=0,\quad
R(\infty)=0.
\end{equation}
\begin{thrm}[Weinstein~\cite{Weinstein-83}] A sufficient condition for global existence in
the critical NLS~(\ref{eq:dDim_CriticalNLS}) is~$\|\psi_0\|_2^2<\Pcr$, where~$\Pcr=\|R^{(0)}\|_2^2$,
and~$R^{(0)}$ is the ground state of
equation~(\ref{eq:dDimCritical_R_ODE}).
\end{thrm}

The critical NLS~(\ref{eq:dDim_CriticalNLS}) admits the explicit
solution
\begin{subequations}
 \label{eq:ExplicitBlowupSolution}
\begin{equation}\label{eq:ExplicitBlowUpSolutions_1}
\psiex(t,r)=\frac1{L^{d/2}(t)}R^{(0)}\left(\frac{r}{L(t)}\right)e^{i\tau+i\frac{L_t}{L}\frac{r^2}{4}},
\end{equation}
where
\begin{equation}\label{eq:ExplicitBlowupSolution_2}
L(t)=T_c-t,\qquad \tau(t)=\int_0^t\frac1{L^2(s)} \,
ds=\frac1{T_c-t}.
\end{equation}
\end{subequations}
More generally, applying the dilation transformation with $\lambda =
\alpha$ and the temporal translation $T_c\longrightarrow \alpha^2
T_c$ shows that the critical NLS~(\ref{eq:dDim_CriticalNLS}) admits
the explicit solutions
\begin{subequations}\label{eq:ExplicitBlowupSolutionAlpha}
\begin{equation}\label{eq:ExplicitBlowupSolutionAlpha_1}
\psiexalpha(t,r)=\frac1{L_\alpha^{d/2}(t)}R^{(0)}\left(\frac{r}{L_\alpha(t)}\right)e^{i\tau_\alpha+i\frac{(L_\alpha)_t}{L_\alpha}\frac{r^2}{4}},
\end{equation}
where
\begin{equation}\label{eq:ExplicitBlowupSolutionAlpha_2}
L_\alpha(t)=\alpha(T_c-t),\qquad
\tau_\alpha(t)=\int_0^t\frac1{L_\alpha^2(s)} \,
ds=\frac1{\alpha^2}\frac1{T_c-t},\qquad\alpha>0.
\end{equation}
\end{subequations}
The explicit
solutions~(\ref{eq:ExplicitBlowupSolution})--(\ref{eq:ExplicitBlowupSolutionAlpha})
become singular at~$t=T_c$. These solutions are unstable, however,
as the have exactly the critical power for collapse. Therefore, any
infinitesimal perturbation which decreases their power, will arrest
the collapse.

When a solution of the critical NLS, whose power is slightly
above~$\Pcr$, undergoes a stable collapse, it splits into two
components: A collapsing core that approaches the universal
$\psi_R$~profile and blows up at the loglog law rate, and a
non-collapsing tail~($\phi$) that does not participate in the
collapse process:
\begin{thrm}[Merle and Raphael~\cite{Merle-Raphael-2003}, \cite{Merle-04}, \cite{Merle-05}, \cite{Merle-05b},
\cite{Merle-06}, \cite{Merle-06b},\cite{Raphael-05}]
\label{thm:Merle+Raphael}
 Let $d=1,2,3,4,5$, and let~$\psi$ be a
solution of the critical NLS~(\ref{eq:dDim_CriticalNLS}) that
becomes singular at~$T_c$. Then, there exists a universal
constant~$m^*>0$, which depends only on the dimension, such that for
any $\psi_0 \in H^1$ such that
$$
\Pcr <\||\psi_0||_2^2< \Pcr+m^*, \qquad
   H_G(\psi_0):= H(\psi_0)-\left ( \frac{\mbox{Im} \int \psi_0^* \nabla \psi_0  }{||\psi_0||_2} \right)^2<0,
$$
the following hold:
\begin{enumerate}
  \item    There exist parameters
  $(\tau(t), {\bf x}_0(t), L(t)) \in {\Bbb R} \times{\Bbb R}^d\times {\Bbb R}^+$,
  and a function $0\neq\phi \in L^2$,  such that
$$
 \psi(t,{\bf x}) -  \psi_R(t,{\bf x}-  {\bf x}_0(t))
 \stackrel{L^2}{\longrightarrow} \phi ({\bf x}), \qquad t \longrightarrow T_c,
$$
where
\begin{equation}\label{eq:PsiR0_Def}
\psi_R(t,{\bf x})  = \frac{1}{L^{d/2}(t)}  R^{(0)} \left(\frac{|{\bf
x}|}{L(t)} \right)e^{i \tau(t)},
\end{equation}
and $R^{(0)}$ is the ground state of
equation~(\ref{eq:dDimCritical_R_ODE}).
\item  As $t \longrightarrow T_c$,
\begin{equation}
\label{eq:loglog_MR} L(t)  \sim \sqrt{2 \pi}
\left(\frac{T_c-t}{\log|\log(T_c-t)|}\right)^{1/2} \!\!\!\!\!\!\!\!
\qquad\mbox{\bf (loglog law)}.
\end{equation}
\end{enumerate}
\end{thrm}

\subsection{Supercritical NLS}\label{subsec:ReviewSupercriticalNLS}
In contrast to the extensive theory on singularity formation in the
critical NLS, much less is known about the supercritical case
\begin{equation}\label{eq:1D_SC_NLS}
    i\psi_t(t,\x)+\Delta\psi+|\psi|^{p-1}\psi=0,  \qquad (p-1)d>4.
\end{equation}
Numerical simulations and formal calculations (see~\cite[Chapter
7]{Sulem-99} and the references therein), and recent rigorous
analysis in the slightly-supercritical
regime~$0<\frac{(p-1)d}{2}-2\ll1$~\cite{Merle-10} show that
peak-type singular solutions of the supercritical
NLS~(\ref{eq:1D_SC_NLS}) collapse with a self-similar asymptotic
profile~$\psi_Q$, where
\begin{equation}\label{eq:psiS}
\psi_Q(
t,r)=\frac{1}{L^{2/(p-1)}(t)}Q\left(\rho\right)e^{i\tau+i\frac{L_t}{L}r^2},
\qquad\rho=\frac{r}{L(t)},\qquad \tau=\int_0^t\frac{ds}{L^2(s)}.
\end{equation}
The blowup rate of~$L(t)$ is a square root, i.e.,
\begin{equation}\label{eq:L_SuperCritical}
L( t)\sim \LC\sqrt{ T_c- t},\qquad  t\to T_c,
\end{equation} where $\LC>0$. In addition, the self-similar
profile~$Q$ is the zero-Hamiltonian, monotonically-decreasing
solution of
\begin{equation}\label{eq:Q_ODE}
Q''(\rho)-\left(1+i\frac{p-5}{4(p-1)}\kappa^2-\frac{\kappa^4}{16}\rho^2\right)Q+|Q|^{p-1}Q=0,\qquad
Q'(0)=0.
\end{equation}
\section{Effect of linear and nonlinear damping - review}\label{sec:NL_DampingEffect}

In~\cite{Fibich-2001}, Fibich studied asymptotically and numerically
the effect of damping on blowup in the critical NLS. He showed that
when the damping is linear, i.e.,
\begin{equation}\label{eq:CriticalLinearDampedLS}
i\psi_t(t,\textbf{x})+\Delta\psi+|\psi|^{4/d}\psi+i\delta\psi=0,\qquad\psi(0,\x)=\psi_0(\x),
\end{equation}
if the initial condition~$\psi_0(\x)$ is such that the solution
of~(\ref{eq:CriticalLinearDampedLS}) becomes singular
for~$\delta=0$, then the solution
of~(\ref{eq:CriticalLinearDampedLS}) exists globally only
if~$\delta$ is above a threshold value~$\delta_c>0$ (which depends
on~$\psi_0$). Therefore, {\em{linear damping cannot play the role of
"viscosity" in continuations of solutions of the NLS}}. When,
however, the damping exponent is critical or supercritical, i.e.,
\begin{equation}\label{eq:CriticalDampedNLS}
i\psi_t(t,\textbf{x})+\Delta\psi+\left(1+i\delta\right)|\psi|^{4/d}\psi=0,\qquad0<\delta\ll1,
\end{equation}
or
\begin{equation}\label{eq:SuperCriticalDampedNLS}
i\psi_t(t,\textbf{x})+\Delta\psi+|\psi|^{4/d}\psi+i\delta|\psi|^{q-1}\psi=0,\qquad0<\delta\ll1\qquad
q-1>4/d,
\end{equation}
respectively, then regardless of how small~$\delta$ is, collapse is
always arrested. Therefore, Fibich suggested that nonlinear damping
can "play the role of viscosity" in defining weak NLS solutions,
i.e., we can define the continuation
\begin{equation}\label{eq:DampedNLS_Continuation}
\psi:=\lim_{\delta\rightarrow 0+}\psi^{(\delta)},
\end{equation}
where~$\psi^{(\delta)}$ is the solution
of~(\ref{eq:CriticalDampedNLS})
or~(\ref{eq:SuperCriticalDampedNLS}).

Passot, Sulem and Sulem proved that high-order nonlinear damping
always prevents collapse for~$d=2$. Antonelli and Sparber extended
this result to~$d=1$ and~$d=3$:
\begin{thrm}[\cite{Sulem-Passot-2005,Antoneli-Sparber-2010}]
\label{thrm:SulemThrm}
 The d-dimensional
cubic NLS with nonlinear damping
\begin{equation}\label{eq:CauchyProblemNLSCubic_1}
i\psi_t+\Delta\psi+\lambda|\psi|^2\psi+i\delta|\psi|^{q-1}\psi=0,\qquad\lambda\in\mathbb{R},\quad\delta>
0,
\end{equation}
where~$\psi_0(\x)\in H^1(\mathbb{R}^d)$,  $3<q<\infty$ if~$d=1,2$,
and~$3<q<5$ if~$d=3$,
has a unique global in-time solution.
\end{thrm}
This rigorously shows that high-order nonlinear damping can play the
role of "viscosity". More recently, Antonelli and Sparber proved
global existence for the case where the damping exponent is equal to
that of the nonlinearity:
\begin{thrm}[\cite{Antoneli-Sparber-2010}]
\label{thrm:AntonelliThrm}
 Consider the cubic nonlinear NLS with
a cubic nonlinear damping
\begin{equation}\label{eq:CauchyProblemNLSCubic}
i\psi_t(t,\x)+\Delta\psi+(1+i\delta)|\psi|^2\psi=0,
\end{equation}
where~$\psi_0(\x)\in H^1(\mathbb{R}^d),~\x\psi_0\in
L^2(\mathbb{R}^d)$, and~$d\leq3$. Then, for any~$\delta\geq1$,
equation~(\ref{eq:CauchyProblemNLSCubic}) has a unique global
in-time solution.
\end{thrm}
Theorem~\ref{thrm:AntonelliThrm} does not show that critical
nonlinear damping can play the role of viscosity. We note, however,
that the asymptotic analysis and simulations
of~\cite{Fibich-Klein-2010,Fibich-2001} strongly suggest that the
solution of~(\ref{eq:CriticalDampedNLS}) exists globally for
any~$0<\delta\ll1$.

\subsection{Explicit continuation
of~$\psiex$}\label{subsec:ExplicitContDampedNLS}

In~\cite{Fibich-Klein-2010}, Fibich and Klein calculated explicitly
the vanishing nonlinear-damping
limit~(\ref{eq:DampedNLS_Continuation}) of the explicit
solution~$\psi_{\rm{explicit}}$:
\begin{prop}[\cite{Fibich-Klein-2010}]
\label{Proposition:DampedNLSWeakSol}
Let~$\psi^{(\delta)}(t,r)$ be the solution
of the NLS~(\ref{eq:CriticalDampedNLS}) with the initial condition
\begin{equation}\label{eq:LemmaIC}
\psi_0(r)=\psiex(0,r),
\end{equation}
see~\eqref{eq:ExplicitBlowupSolution}. Then, for
any~$\theta\in\mathbb{R}$, there exists a
sequence~$\delta_n\rightarrow0+$ (depending on~$\theta$), such that
\begin{equation}\label{eq:DampedWeakSolution}
\lim_{\delta_n\rightarrow 0+}\psi^{(\delta_n)}(t,r)=\left\{
\begin{array}{l l}
  \psiex(t,r) & \quad~0\leq t<T_c,\\
  \psi_{\rm explicit, \kappa}^*(2T_c-t,r)e^{i\theta} & \quad T_c<t<\infty,
\end{array} \right.
\end{equation}
where~$\psi_{\rm explicit, \kappa}$ is given
by~(\ref{eq:ExplicitBlowupSolutionAlpha_1}) with $\alpha = \kappa$,
\begin{equation}\label{eq:AlphaExpression}
\kappa=\pi\left[B_i(0)A'_i(s^*)-A_i(0)B'_i(s^*)\right]\approx 1.614,
\end{equation}
$A_i(s)$ and~$B_i(s)$ are the Airy and Bairy functions,
respectively, and~$s^*\approx-2.6663$ is the first negative root
of~$G(s)=\sqrt{3}A_i(s)-B_i(s)$.

In particular, the limiting width of the solution is given by
\begin{equation}\label{eq:DampedNLS_L_weak}
\lim_{\delta\rightarrow0+}L(t)=\left\{
\begin{array}{l l}
  T_c-t & \quad 0\leq t<T_c,\\
  \kappa(t-T_c) & \quad T_c<t<\infty,
\end{array} \right.
\end{equation}
\end{prop}
Therefore, the continuation is also an explicit
minimal-power solution, but with a higher (defocusing) velocity. In
addition, the solution beyond~$T_c$ is only determined up to 
multiplication by an unknown constant phase~$e^{i\theta}$.

\subsection{Continuation for loglog collapse}

In~\cite{Fibich-Klein-2010}, Fibich and Klein showed asymptotically
and numerically that the continuation of solutions that undergo a
loglog collapse has an infinite-velocity expanding core, that is
determined up to a multiplicative constant phase~$e^{i\theta}$:
\begin{prop}[\cite{Fibich-Klein-2010}]
\label{conject:DampedNLS_Conject}
 Let $\psi_0(r)$ be a radial
initial condition, such that the corresponding solution~$\psi$ of
the undamped critical NLS~(\ref{eq:dDim_CriticalNLS}) collapses with
the $\psi_R$~profile at the loglog law blowup rate at~$T_c$. Let
$\psi^{(\delta)}$ be the solution of the damped
NLS~(\ref{eq:CriticalDampedNLS}) with the same initial condition.
Then,
$$
\lim_{\delta \to 0+} \psi^{(\delta)} = \psi, \qquad 0 \le t < T_c.
$$
In addition, for any~$0<\delta\ll1$, there
exists~$\theta(\delta)\in\mathbb{R}$, and a function~$\phi\in L^2$,
such that
$$
\lim_{\delta\rightarrow
0+}\left[\psi^{(\delta)}(t,r)-{\rev{\psi^\ast_{R}}}(2T_c-t,r;\delta)e^{i\theta(\delta)}\right]\stackrel{L^2}\longrightarrow\phi(r),
\qquad t\longrightarrow T_c+,
$$
{\rev{where~$\psi_R$ is given by~(\ref{eq:PsiR0_Def}) with some
function~$L(t;\delta)$, such that}}
\[
{\rev{\lim_{t\rightarrow
T_c+}\lim_{\delta\rightarrow0+}L(t;\delta)=0,\qquad\lim_{t\rightarrow
T_c+}\lim_{\delta\rightarrow0+}L_t(t;\delta)=\infty}},
\qquad\lim_{\delta\rightarrow0+}\theta(\delta)=\infty.
\]
\end{prop}

\begin{rmk}
We use the terminology {\em Continuation Result}, in order to emphasize that 
the proofs of Continuation Results~\ref{Proposition:DampedNLSWeakSol}
and~\ref{conject:DampedNLS_Conject} are based on asymptotic analysis
and numerical simulations, and are not rigorous.
\end{rmk}

\section{The critical damping
exponent}\label{sec:CriticalDampingExponent}

The rigorous theorems in Section~\ref{sec:NL_DampingEffect} suggest
that solutions of the NLS~(\ref{eq:DampedNLS}) with~$\psi_0(\X)\in
H^1$ always exist globally when~$p<q$. These theorems, however, do
not cover the case~$p=q$. The asymptotic analysis and simulations
of~\cite{Fibich-Klein-2010,Fibich-2001} strongly suggest that in the
critical case, when~$p=q$ the solution always exists globally, see
Continuation Results~\ref{Proposition:DampedNLSWeakSol}
and~\ref{conject:DampedNLS_Conject}. Since there are no rigorous and
asymptotic results for the supercritical case with~$p=q$, we study
this case numerically.

\begin{figure}[ht!]
\begin{center}
\scalebox{0.7}{\includegraphics{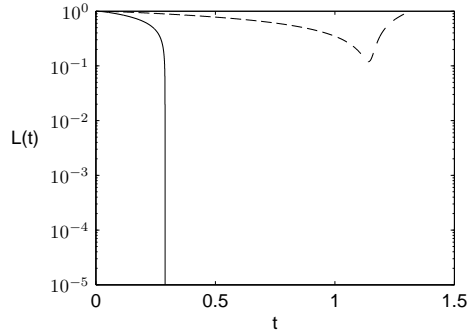}}
\caption{Solution of~(\ref{eq:DampedNLS_EqualExponents})
for~$\delta=5\cdot10^{-3}$ with~$p=5$ (dashes) and~$p=7$
(solid).}\label{fig:Critical_vs_Supercritical_DampedNLS}
\end{center}
\end{figure}

Consider the one-dimensional NLS with~$p=q$
\begin{subequations}\label{eq:DampedNLS_EqualExponents}
\begin{equation}\label{eq:DampedNLS_EqualExponents_1}
i\psi_t(t,x)+\psi_{xx}+(1+i\delta)|\psi|^{p-1}\psi=0,
\end{equation}
with the perturbed solitary-wave initial condition
\begin{equation}\label{eq:DampedNLS_EqualExponents_2}
\psi_0(x)=1.05R_p(x),
\end{equation}
\end{subequations}
where~$R_p(x)$ is given by~(\ref{eq:1D_ground_state_solutions}). In
Figure~\ref{fig:Critical_vs_Supercritical_DampedNLS} we
solve~(\ref{eq:DampedNLS_EqualExponents})
with~$\delta=5\cdot10^{-3}$, and plot
\begin{equation}
   \label{eq:BeamWidthDefinitions_1}
L(t)=\left|\frac{\psi(0,0)}{\psi(t,0)}\right|^{(p-1)/2}.
\end{equation}
In the critical case~$p=5$, the collapse is arrested after focusing
by~$\approx10$. In the supercritical case~$p=7$, however, the
collapse is not arrested after focusing by~$10^{5}$. This and
similar simulations suggest that unlike the critical case, {\em{in
the supercritical case, the condition~$p<q$ is necessary for
ensuring global existence in~(\ref{eq:DampedNLS}).}}

\section{Supercritical NLS}
\label{sec:SuperCriticalNLS}

We now consider the effect of small nonlinear damping in the
supercritical NLS. Let~$\psi^{(\delta)}(t,x)$ be the solution of the
one-dimensional supercritical damped NLS ($d=1$,~$p=7$,~$q=9$)
\begin{subequations}\label{eq:SC_DampedNLS_GaussianIC}
\begin{equation}\label{eq:SC_DampedNLS__GaussianIC_1}
i\psi_t(t,x)+\psi_{xx}+|\psi|^6\psi+i\delta|\psi|^8\psi=0,
\end{equation}
with the initial condition
\begin{equation}\label{eq:SC_DampedNLS_GaussianIC_2}
\psi_0(x)=1.3e^{-x^2}.
\end{equation}
\end{subequations}
Let
\begin{equation}\label{eq:T_delta_max_def}
T^{(\delta)}_{\max}=\argmax_t\|\psi^{(\delta)}(t,x)\|_\infty,
\end{equation}
denote the time at which the focusing is maximal.

\begin{figure}[ht!]
\begin{center}
\scalebox{0.7}{\includegraphics{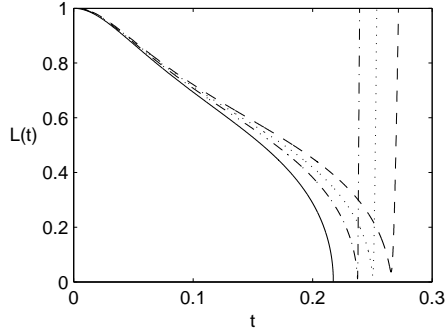}}
\caption{Solution of~(\ref{eq:SC_DampedNLS_GaussianIC})
for~$\delta=0$ (solid),~$\delta=5\cdot10^{-3}$
(dashes-dots),~$\delta=7.5\cdot10^{-3}$ (dots), and~$\delta=10^{-2}$
(dashes).}\label{fig:SC_DampedNLS_BeamWidth_GaussianIC}
\end{center}
\end{figure}

In Figure~\ref{fig:SC_DampedNLS_BeamWidth_GaussianIC} we plot the
solution of~(\ref{eq:SC_DampedNLS_GaussianIC}) for various values
of~$\delta$. In all cases, the collapse is arrested in a highly
asymmetric way with respect to~$T_{\rm arrest}^{(\delta)}$. In addition,
the post-collapse defocusing rate of the singular core "appears" to
increase to infinity as~$\delta\rightarrow0+$.

\begin{figure}[ht!]
\begin{center}
\scalebox{0.8}{\includegraphics{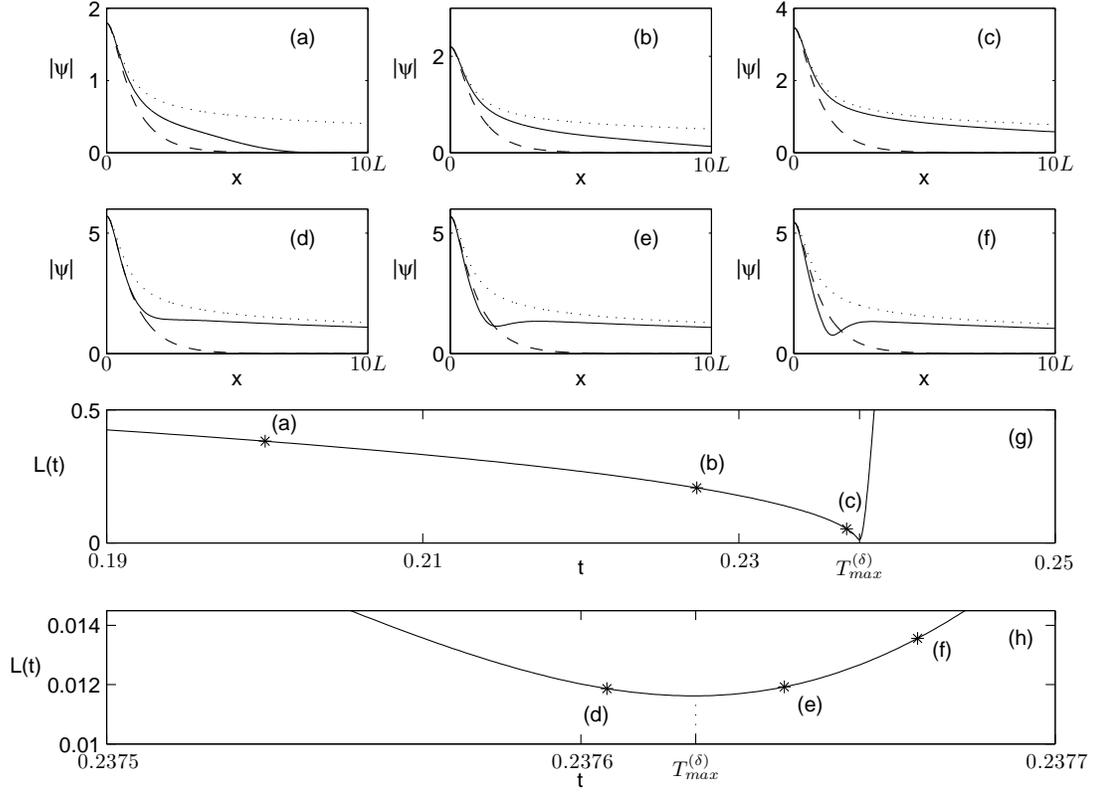}}
\caption{Solution of~(\ref{eq:SC_DampedNLS_GaussianIC})
for~$\delta=5\cdot10^{-3}$ (solid), and the fitted~$|\psi_Q|$ (dots)
and~$|\psi_R|$ (dashes).
(a)~$t\approx0.2$,~$L\approx0.38$~(b)~$t\approx0.227$,~$L\approx0.2$
~(c)~$t\approx0.2368$,~$L\approx0.05$~(d)~$t\approx0.2376$,~$L\approx0.0118$
~(e)~$t\approx0.23764$,~$L\approx0.0119$
~(f)~$t\approx0.23767$,~$L\approx0.135$. The stars in~(g) and~(h)
denote the values of~$t$ and~$L(t)$ for the data in
subplots~(a)--(f).}\label{fig:ProfileFittingDampedSCNLS_GaussianIC}
\end{center}
\end{figure}

In Figure~\ref{fig:ProfileFittingDampedSCNLS_GaussianIC} we compare
the profile of the solution of~(\ref{eq:SC_DampedNLS_GaussianIC})
with~$\delta=10^{-3}$ with the supercritical~$\psi_Q$
and~$\psi_R$~profiles, where
\begin{subequations}\label{eq:AsymptoticProfiles}
\begin{equation}\label{eq:PsiR_Def}
|\psi_R(t,x)|=\frac1{L_R^{2/(p-1)}(t)}R\left(\frac{x}{L_R(t)}\right),\qquad
L_R(t)=\left|\frac{R(0)}{\psi^{(\delta)}(t,0)}\right|^{(p-1)/2},
\end{equation}
\begin{equation}
|\psi_Q(t,x)|=\frac1{L_Q^{2/(p-1)}(t)}\left|Q\left(\frac{x}{L_Q(t)}\right)\right|,\qquad
L_Q(t)=\left|\frac{Q(0)}{\psi^{(\delta)}(t,0)}\right|^{(p-1)/2},
\end{equation}
\end{subequations}
and R and Q are the solutions of~(\ref{eq:dDim_R_ODE})
and~(\ref{eq:Q_ODE}), respectively, with~$d=1$ and~$p=7$.
The NLS solution initially approaches the~$\psi_Q$ profile, see
Figure~\ref{fig:ProfileFittingDampedSCNLS_GaussianIC}(a--c). This is
to be expected, since when~$\delta=0$ the solution collapses with
the~$\psi_Q$ profile, see
Section~\ref{subsec:ReviewSupercriticalNLS}. As the solution
approaches~$T_{\rm arrest}^{(\delta)}$, however, the collapsing core moves
away from~$\psi_Q$ and towards~$\psi_R$, see
Figure~\ref{fig:ProfileFittingDampedSCNLS_GaussianIC}(d), and it
remains close to~$\psi_R$ for a "short time"
after~$T_{\rm arrest}^{(\delta)}$, see
Figure~\ref{fig:ProfileFittingDampedSCNLS_GaussianIC}(e).
Eventually, as the collapsing core continues to defocus, it
interacts with its tail and "loses" its~$\psi_R$ profile, see
Figure~\ref{fig:ProfileFittingDampedSCNLS_GaussianIC}(f).

\begin{figure}[ht!]
\begin{center}
\scalebox{0.7}{\includegraphics{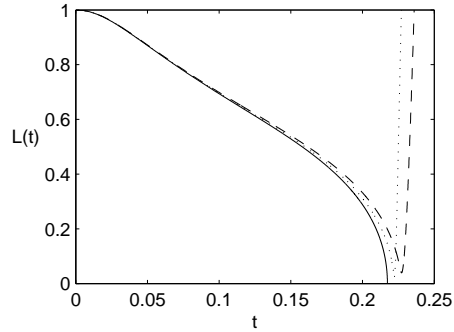}}
\caption{Solution of~(\ref{eq:SC_DampedNLS_DampedOrd11})
for~$\delta=0$ (solid),~$\delta=5\cdot10^{-4}$ (dots),
and~$\delta=10^{-3}$
(dashes).}\label{fig:SC_DampedNLS_BeamWidth_GaussianIC_DampedOrd11}
\end{center}
\end{figure}

\begin{figure}[ht!]
\begin{center}
\scalebox{0.8}{\includegraphics{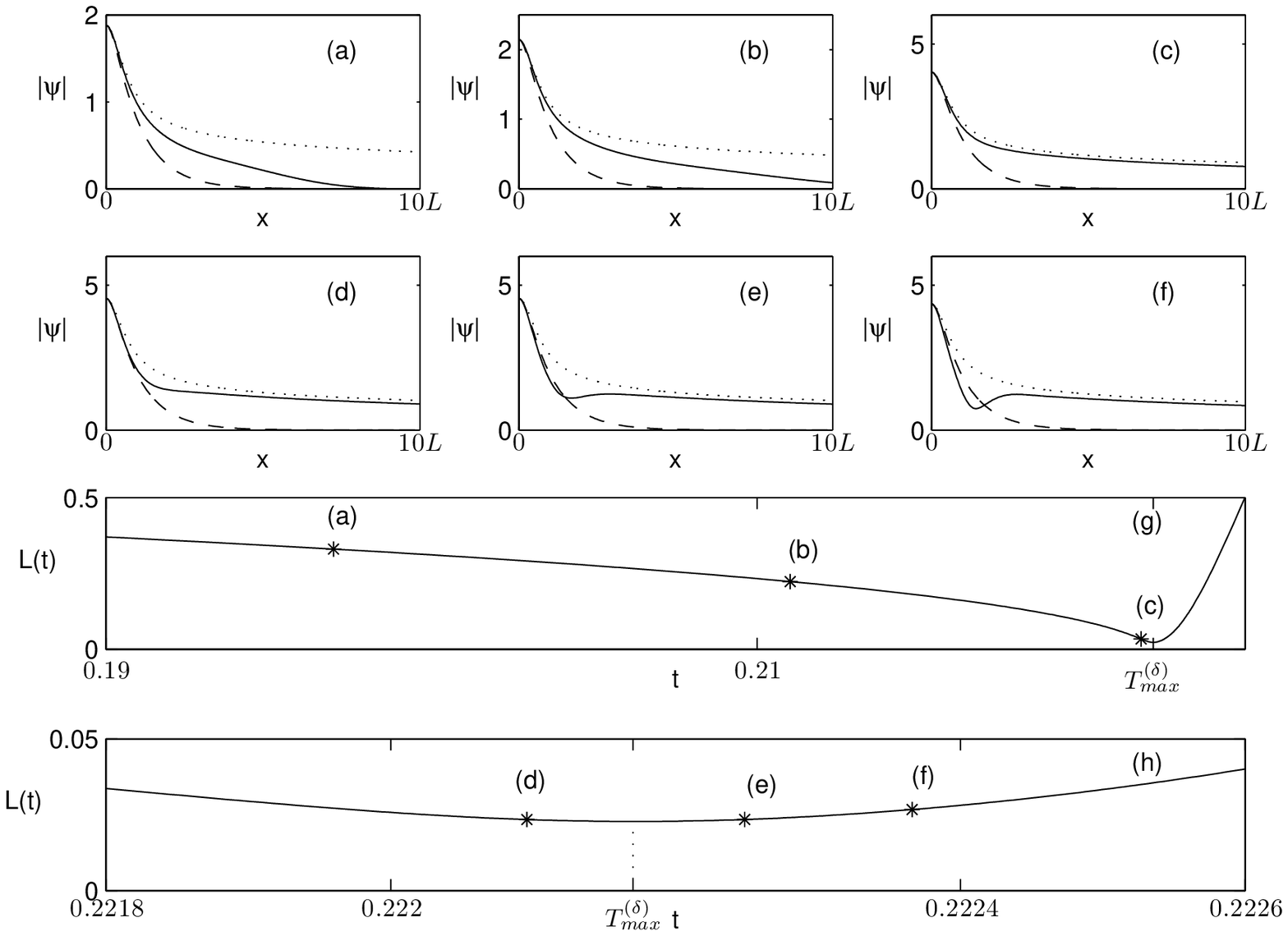}}
\caption{Same as
Figure~\ref{fig:ProfileFittingDampedSCNLS_GaussianIC} for the
solution of~(\ref{eq:SC_DampedNLS_DampedOrd11}).
(a)~$t\approx0.19$,~$L\approx0.33$~(b)~$t\approx0.211$,~$L\approx0.22$
~(c)~$t\approx0.221$,~$L\approx0.033$~(d)~$t\approx0.222$,~$L\approx0.02341$
~(e)~$t\approx0.2222$,~$L\approx0.02343$
~(f)~$t\approx0.2223$,~$L\approx0.0267$.}\label{fig:ProfileFittingDampedSCNLS_GaussianIC_DampedOrd11}
\end{center}
\end{figure}

Next, we repeat the above simulation with a higher nonlinear damping
exponent ($q=11$). Specifically, we solve the NLS
\begin{subequations}\label{eq:SC_DampedNLS_DampedOrd11}
\begin{equation}\label{eq:SC_DampedNLS_DampedOrd11_1}
i\psi_t(t,x)+\psi_{xx}+|\psi|^6\psi+i\delta|\psi|^{10}\psi=0,
\end{equation}
with the initial condition
\begin{equation}\label{eq:SC_DampedNLS_DampedOrd11_2}
\psi_0(x)=1.3e^{-x^2}.
\end{equation}
\end{subequations}
Figure~\ref{fig:SC_DampedNLS_BeamWidth_GaussianIC_DampedOrd11}
and~\ref{fig:ProfileFittingDampedSCNLS_GaussianIC_DampedOrd11} show
that the qualitative behavior of the solution is exactly the same as
that of the solution of~(\ref{eq:SC_DampedNLS_GaussianIC}).

Therefore, we conclude that solutions of the supercritical
NLS~(\ref{eq:DampedNLS}) with~$q>p>1+4/d$ and~$0<\delta\ll1$:
\begin{enumerate}
\item Exist globally.
\item Are highly asymmetric with respect to~$T_{\rm arrest}^{(\delta)}$.
\item The post-collapse velocity of the defocusing core goes to infinity as~$\delta\rightarrow0+$.
\item The asymptotic profile around~$T_{\rm arrest}^{(\delta)}$ is~$\psi_R$, and
not~$\psi_Q$. 
   
    The fact that as~$t\rightarrow T_{\rm arrest}^{(\delta)}$ the
profile is not given by~$\psi_Q$ is not surprising, since the
nonlinear damping perturbation obviously has a significant effect
near~$T_{\rm arrest}^{(\delta)}$, and therefore there is no reason why it
should not change the solution profile. What is surprising is that
the profile changes to the supercritical~$\psi_R$ profile. As far as
we know, this is the first observation in which the asymptotic
profile in the supercritical NLS is given by the supercritical $\psi_R$~profile.\footnote{The standing-ring solutions of the undamped supercritical  NLS with~$p=5$ and~$d>1$ also collapse 
with the $\psi_R$~profile~\cite{Raphael-06,SC_rings-07,Raphael-08,Standing-10}. 
In that case, however, $\psi_R$ is the asymptotic profile of the {\em critical} one-dimensional quintic NLS.}

\end{enumerate}
\section{Critical NLS}\label{sec:CriticalNLS}
\subsection{Continuation of loglog collapse}
In~\cite{Fibich-Klein-2010}, we studied the effect of nonlinear
damping in the critical NLS with~$p=q$ with initial conditions that
lead to a loglog collapse, see
Continuation Result~\ref{conject:DampedNLS_Conject}. We now consider the
case~$p>q$.

Consider the damped one-dimensional critical NLS
($d=1$,~$p=5$,~$q=7$)
\begin{subequations}\label{eq:CriticalNLS_SC_Damping}
\begin{equation}\label{eq:CriticalNLS_SC_Damping_1}
i\psi_t(t,x)+\psi_{xx}+|\psi|^4\psi+i\delta|\psi|^6\psi=0,
\end{equation}
with the initial condition
\begin{equation}\label{eq:CriticalNLS_SC_Damping_2}
\psi_0(x)=1.6e^{-x^2},
\end{equation}
\end{subequations}
whose power is~$4\%$ above the critical power for collapse.
When~$\delta=0$, the NLS solution collapses with the~$\psi_R$
profile at the loglog blowup rate.

\begin{figure}[ht!]
\begin{center}
\scalebox{0.8}{\includegraphics{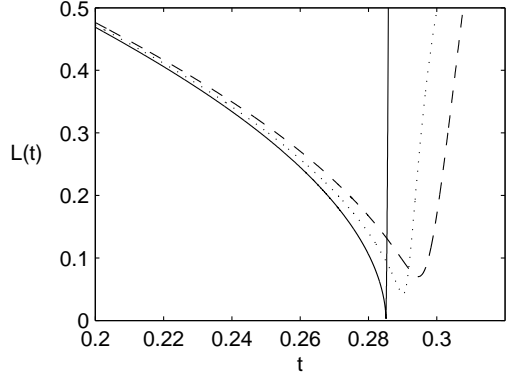}}
\caption{Solution of~(\ref{eq:CriticalNLS_SC_Damping})
for~$\delta=10^{-5}$ (solid),~$\delta=2.5\cdot10^{-4}$ (dots) ,
and~$\delta=5\cdot10^{-4}$
(dashes).}\label{fig:CriticalNLS_SC_Damping_L}
\end{center}
\end{figure}

In Figure~\ref{fig:CriticalNLS_SC_Damping_L} we
solve~(\ref{eq:CriticalNLS_SC_Damping}) for various values
of~$\delta$. In all cases, the collapse is arrested in a highly
asymmetric way with respect to~$T_{\rm arrest}^{(\delta)}$. In addition,
the post-collapse defocusing rate appears to increase to infinity
as~$\delta\rightarrow0+$. This qualitative behavior is as in the
case~$p=q$, see Continuation Result~\ref{conject:DampedNLS_Conject}. Therefore, we
conclude that the qualitative behavior for~$q=p$ and for~$q>p$ is
the same.

\begin{figure}[ht!]
\begin{center}
\scalebox{0.8}{\includegraphics{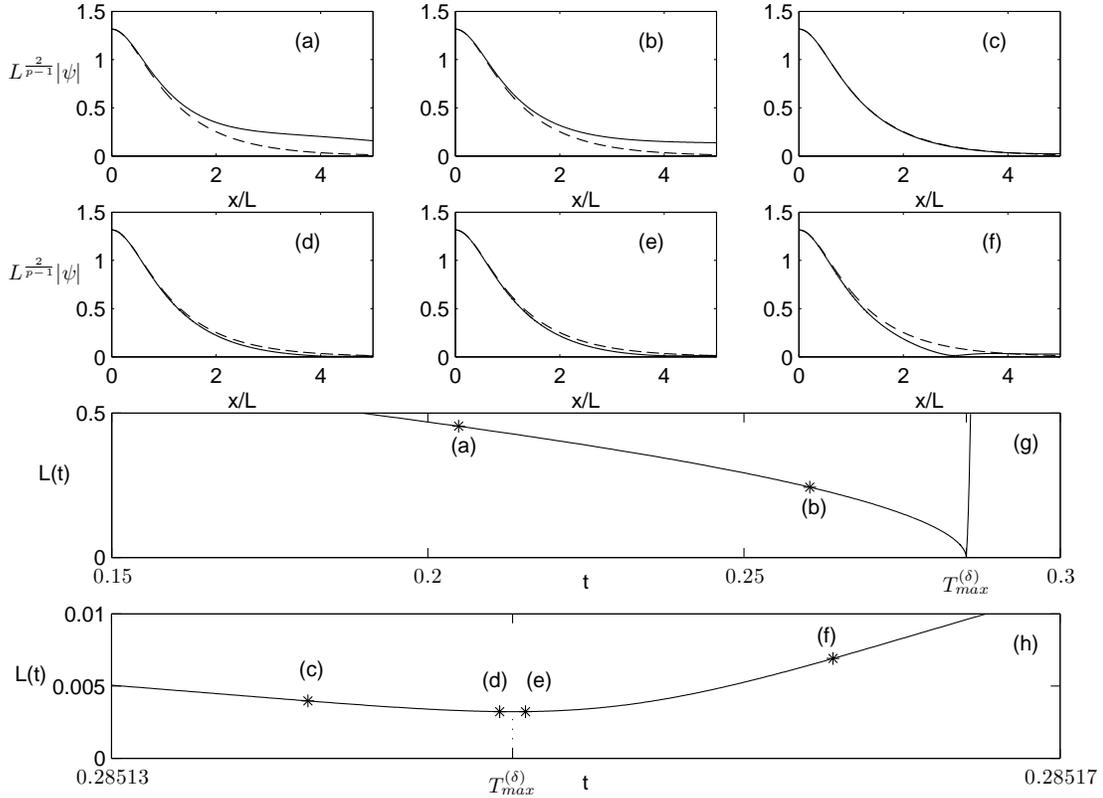}}
\caption{Solution of~(\ref{eq:CriticalNLS_SC_Damping})
for~$\delta=10^{-5}$ (solid), and the fitted~$|\psi_R|$ (dashes).
(a)~$t\approx0.205$,~$L\approx0.45$~(b)~$t\approx0.26$,~$L\approx0.24$
~(c)~$t\approx0.28514$,~$L\approx0.0037$~(d)~$t\approx0.285146$,~$L\approx0.0032$
~(e)~$t\approx0.285147$,~$L\approx0.0032$
~(f)~$t\approx0.28515$,~$L\approx0.0046$. The stars in~(g) and~(h)
denote the values of~$t$ and~$L(t)$ for the data in
subplots~(a)--(f).}\label{fig:ProfileFittingDampedCriticalNLS_GaussianIC_DampedOrd7}
\end{center}
\end{figure}

In
Figure~\ref{fig:ProfileFittingDampedCriticalNLS_GaussianIC_DampedOrd7}
we compare the profile of the solution
of~(\ref{eq:CriticalNLS_SC_Damping}) with~$\delta=10^{-5}$ with the
best-fitting critical~$\psi_R$ profile, see~(\ref{eq:PsiR_Def}). The
NLS solution initially approaches the~$\psi_R$ profile, see
Figure~\ref{fig:ProfileFittingDampedCriticalNLS_GaussianIC_DampedOrd7}(a--c).
This is to be expected, since when~$\delta=0$ the solution collapses
with the~$\psi_R$ profile, see Theorem~\ref{thm:Merle+Raphael}. As
the solution approaches~$T_{\rm arrest}^{(\delta)}$, however, the
collapsing core moves away from~$\psi_R$, see
Figure~\ref{fig:ProfileFittingDampedCriticalNLS_GaussianIC_DampedOrd7}(d--e).
Unlike the supercritical case, however, the solution profile 
near~$T_{\rm arrest}^{(\delta)}$ is still ``close'' to~$\psi_R$. This is because in the critical case, perturbations arrest the collapse when they are still small compared with the nonlinearity and diffraction~\cite{PNLS-99}. 
Eventually, as the collapsing core continues to defocus, it
interacts with its tail and "loses" its~$\psi_R$~profile, see
Figure~\ref{fig:ProfileFittingDampedCriticalNLS_GaussianIC_DampedOrd7}(f).

In summary, nonlinearly-damped loglog solutions of
the critical NLS with~$q \ge p$ have the following properties:
\begin{enumerate}
\item The solutions are highly asymmetric with respect to~$T_{\rm arrest}^{(\delta)}$.
\item The post-collapse defocusing velocity goes to infinity as~$\delta\rightarrow0+$.
\item The asymptotic profile near~$T_{\rm arrest}^{(\delta)}$ is "slightly"
different form~$\psi_R$.
\end{enumerate}

\subsection{Continuation of~$\psi_{\rm{explicit}}$}

Consider the critical NLS with nonlinear damping
\begin{subequations}
 \label{eq:1D_GeneralCriticalDampedNLS_IP}
\begin{equation}
 \label{eq:1D_GeneralCriticalDampedNLS}
i\psi_t(t,\X)+\Delta\psi+|\psi|^{4/d}\psi+i\delta|\psi|^{q-1}\psi=0,
 \qquad 0<\delta \ll 1, 
\end{equation}
and the initial condition
\begin{equation}
   \label{eq:PsiExplicitIC}
\psi_0(r)=\psi_{\rm{explicit}}(0,r).
\end{equation}
\end{subequations}
When $\delta = 0$, the solution is given by~$\psi_{\rm{explicit}}$, see
equation~(\ref{eq:ExplicitBlowupSolution}). In~\cite{Fibich-Klein-2010},
we calculated explicitly the continuation of~$\psi_{\rm{explicit}}$ when~$q=p$, 
see Continuation Result~\ref{Proposition:DampedNLSWeakSol}. We now consider the continuation
for~$q \not=p$.\footnote{Since $\psi_{\rm{explicit}}$ has exactly the critical power for collapse, any amount of damping will arrest 
the collapse. Therefore, the continuation of~$\psi_{\rm{explicit}}$ can also be defined for $q<p$.}

As in~\cite{Fibich-Klein-2010}, we can use {\em modulation theory}~\cite{PNLS-99} 
to approximate equation~(\ref{eq:1D_GeneralCriticalDampedNLS_IP}) with 
a reduced system of ordinary-differential equations. 
\begin{lem}
Let 
$$
L(t)=\left|\frac{\psi(0,0)}{\psi(t,0)}\right|^{2/d},
$$
 where~$\psi$ is the solution of equation~(\ref{eq:1D_GeneralCriticalDampedNLS_IP}). Then, 
as~$\delta \longrightarrow 0+$, the evolution of~$L(t)$ is governed by the reduced equations
 \begin{subequations}
  \label{eq:GeneralCriticalDampedNLS_ReducedEq}
\begin{equation}
 \label{eq:GeneralCriticalDampedNLS_ReducedEq_1}
\beta_t(t)=-\frac{2c_q\delta}{M}\frac1{L^{(q-1)d/2}},\qquad
L_{tt}(t)=-\frac{\beta(t)}{L^3},
\end{equation}
subject to the initial conditions
\begin{equation}\label{eq:GeneralCriticalDampedNLS_ReducedEq_3}
\beta(0)=0,\qquad L(0)=0,\qquad L_t(0)=-1,
\end{equation}
\end{subequations}
where
\[
\upsilon(\beta)=\left\{
\begin{array}{l l}
  c_{\nu}e^{-\pi/\sqrt{\beta}}, & \qquad \beta>0,\\
  0, & \qquad \beta\leq0,\\
\end{array} \right.
\]
\[
c_{\nu}=\frac{2A^2_R}{M},\qquad
A_R=\lim_{r\rightarrow\infty}e^rr^{(d-1)/2}R^{(0)}(r),
\]
\[
M=\frac1{4}\int_0^\infty r^2 |R^{(0)}|^2r^{d-1}dr,\qquad
c_q=\|R^{(0)}\|^{q+1}_{q+1},\qquad
\]
and~$R^{(0)}$ is the ground state of~(\ref{eq:dDimCritical_R_ODE}).
\end{lem}
\Beginproof
In~\cite{Fibich-2001} it was shown that the reduced equations for the damped NLS~\eqref{eq:1D_GeneralCriticalDampedNLS}
are given by
\begin{equation}\label{eq:GenDampedReducedEq}
\beta_t(t)=-\frac{\nu(\beta)}{L^2}-\frac{2c_q\delta}{M}\frac1{L^{(q-1)d/2}},\qquad
L_{tt}=-\frac{\beta(t)}{L^3}.
\end{equation}
In addition, the initial
conditions for the reduced equations~(\ref{eq:GenDampedReducedEq})
that correspond to the initial condition~(\ref{eq:PsiExplicitIC})
are~$\beta(0)=0$,~$L(0)=0$, and~$L_t(0)=0$, see~\cite{Fibich-Klein-2010}. Since~$\beta(0)=0$, and
since~$\beta_t<0$, then~$\beta(t)<0$. Hence,~$\nu(\beta)\equiv0$.
Therefore, the reduced equations are given by~\eqref{eq:GeneralCriticalDampedNLS_ReducedEq}.
\Endproof

The reduced equations variable~$L(t)$ is the solution width, and is also inversely proportional to the solution amplitude.
The reduced equations variable~$\beta(t)$ is a measure of the acceleration of~$L(t)$, 
and is also linearly proportional to the excess power above~$\Pcr$ of the collapsing core.

\begin{figure}[ht!]
\begin{center}
\scalebox{0.6}{\includegraphics{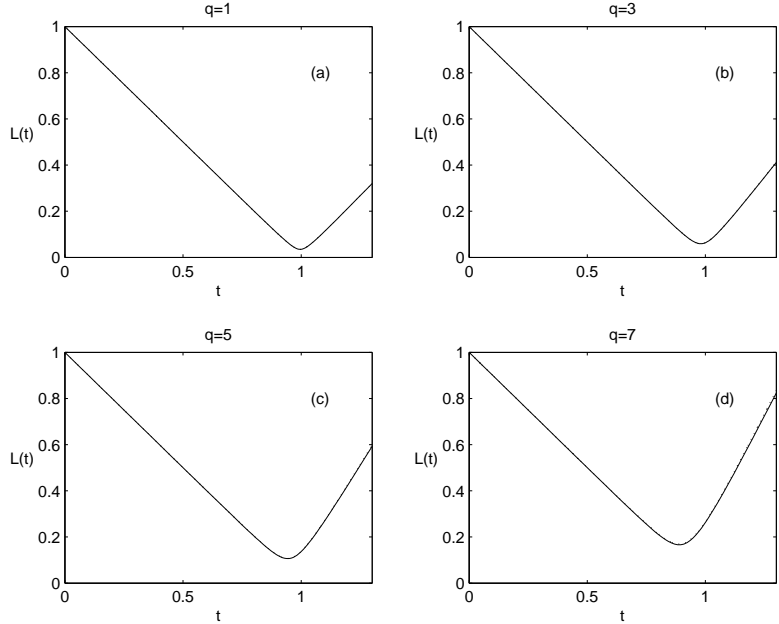}}
\caption{Solution of the reduced
equations~(\ref{eq:GeneralCriticalDampedNLS_ReducedEq}) [solid], and
of the NLS~(\ref{eq:1D_GeneralCriticalDampedNLS_IP}) [dashes],
for~$\delta=2.5\cdot10^{-5}$ and~$d=1$. The two curves are indistinguishable.
(a)~$q=1$.~(b)~$q=3$.~(c)~$q=5$.~(d)~$q=7$.}\label{fig:GeneralCriticalDampedNLS_ReducedEq_vs_Simulation}
\end{center}
\end{figure}

  Since {\em modulation theory} is not rigorous, in Figure~\ref{fig:GeneralCriticalDampedNLS_ReducedEq_vs_Simulation}
we compare the numerical solutions of the reduced
equations~(\ref{eq:GeneralCriticalDampedNLS_ReducedEq}) and the
NLS~(\ref{eq:CriticalNLS_SC_Damping}). This comparison shows that the two solutions are in excellent agreement,
thus providing a strong support to the validity of the reduced equations. 
Therefore, in what follows we study asymptotically and numerically the limit $\delta \longrightarrow 0+$ within the framework of the reduced equations, which is considerably easier than studying the limit $\delta \longrightarrow 0+$  of the nonlinearly-damped NLS.

The extension of Continuation Result~\ref{Proposition:DampedNLSWeakSol} to~$q \not=p$ is as follows.
\begin{prop}
\label{Proposition:DampedNLSWeakSol-new}
Let~$\psi^{(\delta)}(t,r)$ be the solution
of the NLS~(\ref{eq:1D_GeneralCriticalDampedNLS_IP}). Then, for
any~$\theta\in\mathbb{R}$, there exists a
sequence~$\delta_n\rightarrow0+$ (depending on~$\theta$), such that
\begin{equation}\label{eq:DampedWeakSolution}
\lim_{\delta_n\rightarrow 0+}\psi^{(\delta_n)}(t,r)=\left\{
\begin{array}{l l}
  \psiex(t,r) & \quad~0\leq t<T_c,\\
  \psi_{\rm explicit, \kappa(q)}^*(2T_c-t,r)e^{i\theta} & \quad T_c<t<\infty.
\end{array} \right.
\end{equation}
In particular, the limiting width of the solution is given by
\begin{equation}\label{eq:DampedNLS_L_weak}
\lim_{\delta\rightarrow0+}L(t)=\left\{
\begin{array}{l l}
  T_c-t & \quad 0\leq t<T_c,\\
  \kappa(q) (t-T_c) & \quad T_c<t<\infty.
\end{array} \right.
\end{equation}
\end{prop}

\Beginproof We only provide an informal proof, using the reduced equations~\eqref{eq:GeneralCriticalDampedNLS_ReducedEq}. 
 As~$\delta\rightarrow0+$, $\beta_t(t)\rightarrow0$, see equation~\eqref{eq:GeneralCriticalDampedNLS_ReducedEq_1}.
Therefore, since $\beta(0) = 0$, then~$\beta(t)\rightarrow0$. Hence,
$L_{tt}(t)\rightarrow0$. Therefore, $\lim_{\delta\rightarrow0+}L(t)$~is linear in~$t$.
Since~$\lim_{\delta\rightarrow0+}L(T_c)=0$, it follows that
\[
\lim_{\delta\rightarrow0+}L(t;\delta)=\left\{
\begin{array}{l l}
  T_c-t, & \qquad t<T_c,\\
  \kappa(q)(t-T_c), &\qquad T_c<t.\\
\end{array} \right.
\]
The loss of phase follows from the fact that $\lim_{t \to T_c} \arg (\psiex(t,0)) = \infty$  ,see~\cite{Fibich-Klein-2010}.
\Endproof

The result of Continuation Result~\ref{Proposition:DampedNLSWeakSol-new} can be explained as follows.
By continuity, $\lim_{\delta \rightarrow 0+}\psi^{(\delta)}=\psiex$ 
for~$0 \le t<T_c$. The limiting solution 
for~$t>T_c$ is an NLS solution that becomes singular 
 as $t\rightarrow T_c+$, and has exactly the critical power at the singularity. Hence,  
 the limiting solution is a minimal-power solution. Therefore, it has to be given
by~$\psiexalpha$~\cite{Merle-92a,Merle-92b}.

\begin{figure}[ht!]
\begin{center}
\scalebox{0.8}{\includegraphics{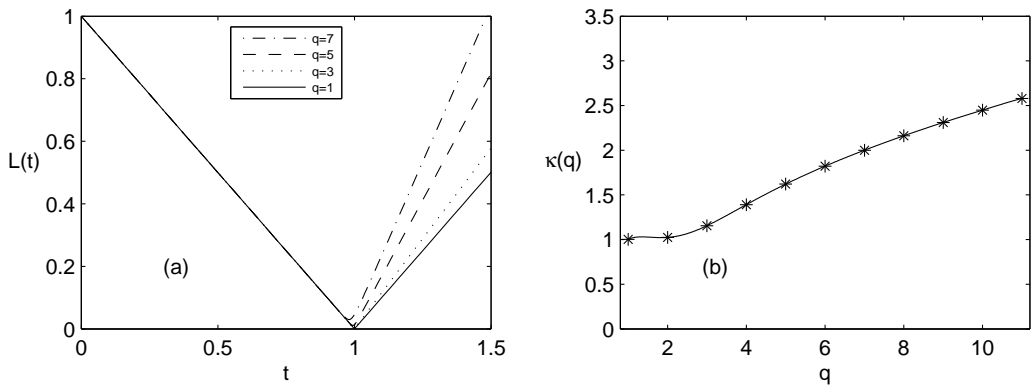}}
\caption{Solution of the reduced
equations~(\ref{eq:GeneralCriticalDampedNLS_ReducedEq})
with~$\delta=10^{-7}$ and various values
of~$q$.~(a)~$L(t)$.~(b)~$\kappa(q)$.}\label{fig:GenDampedCriticalNLS_ReducedEq}
\end{center}
\end{figure}

In Figure~\ref{fig:GenDampedCriticalNLS_ReducedEq} we solve the
reduced equations~(\ref{eq:GeneralCriticalDampedNLS_ReducedEq})
with~$\delta=10^{-7}$, and observe that:
\begin{enumerate}
\item The limiting solutions are indeed linear for~$t<T_c$ and~$t>T_c$.
\item The continuation is asymmetric with respect to~$T_c$.
\item The post-collapse
slope~$\kappa(q)$ increases with~$q$. In~\cite{Fibich-Klein-2010} we
showed that the jump discontinuity in~$\lim_{\delta\rightarrow0+} L_t^2$ at~$T_c$ is related to
the increase of the Hamiltonian as the limiting solution passes
through the singularity. As~$q$ increases, damping affects become
more pronounced, hence there is a larger increase of the
Hamiltonian, hence of the post-collapse slope.
\item When $q=1$, $\kappa(q=1)=1$, i.e., $L(t)$~is symmetric with respect to~$T_c$. 
    Therefore, {\em the linear damping continuation of~$\psiex$ is
symmetric with respect to~$T_c$, even though the problem is not
time-reversible.}

\end{enumerate}

Note that the value of~$\kappa(q=1+4/d)\approx1.614$ was computed
analytically in Continuation Result~\ref{Proposition:DampedNLSWeakSol}.

\section{Final remarks}
  \label{sec:final}

In this study we used numerical simulations to study the effect of 
small nonlinear-damping on singular NLS solutions.  
These simulations suggest that the effect of small nonlinear
damping is qualitatively the same
in the critical NLS with generic initial conditions that lead to a loglog collapse with
the~$\psi_R$ profile, and in the supercritical NLS with generic initial conditions 
that lead to a  square-root collapse with the~$\psi_Q$ profile. 
Moreover, the qualitative effect of nonlinear damping is independent of the value
of~$q$, so long as~$q>p$ in the supercritical case and~$q\geq p$ in
the critical case. Thus, because nonlinear damping destroys the NLS time reversibility, 
the nonlinearly-damped solution is highly
asymmetric with respect to the arrest-of-collapse time~$T_{\rm arrest}^{(\delta)}$. The post-collapse
defocusing velocity~$L_t(t)$ of the singular core goes to infinity as~$\delta \longrightarrow 0+$, 
since the focusing velocity before the singularity goes to infinity for
loglog and square-root blowup rates, and since nonlinear damping
increases the Hamiltonian,\footnote{If we multiply the NLS~\eqref{eq:DampedNLS} by~$\psi_t^*$, 
add the complex-conjugate equation, and integrate by parts, we get that 
$$
H_t = i \delta \int  |\psi|^{q-1} \psi \psi^*_t + c.c.,
$$
where~$c.c.$ stand for complex conjugate. Let $\psi = A e^{i S}$, where~$A$ and~$S$ are real, Then,
$$
H_t = 2 \delta  \int |A|^{q+1} S_t.
$$
 Since for collapsing solutions $S_t \sim L^{-2}(t)$, it follows that $H_t>0$.
}
  hence the ``kinetic energy''.

 Around~$T_{\rm arrest}^{(\delta)}$,
the collapsing core of the singular core  moves away from the
asymptotic profile of the undamped solution.  
 In the supercritical
case the difference between the solution profile
and~$\psi_Q$ for~$t\approx T_{\rm arrest}^{(\delta)}$ is large.
This is intuitive, since damping effects have a large effect when they arrest the collapse.
In the critical case, however, 
the difference between the solution profile and~$\psi_R$
for~$t\approx T_{\rm arrest}^{(\delta)}$ is minor.
This is because critical collapse has the unique property that it can be arrested 
by small perturbations~\cite{PNLS-99}.

{\em Surprisingly,  in the supercritical
case the profile of the nonlinearly-damped solution near~$T_{\rm arrest}^{(\delta)}$ appears to be
given by the supercritical~$\psi_R$ profile}. To the best of our
knowledge, this is the first observation of a solution of the
supercritical NLS that approaches the
supercritical~$\psi_R$~profile.

\subsubsection*{Acknowledgment}
 This research was partially supported by grant~$1023/08$ from the Israel Science
Foundation (ISF).


\end{document}